\RequirePackage{fix-cm}
\documentclass[onecolumn]{svjour3questa}       
\usepackage{amsmath,amssymb}
\usepackage{mathptmx}      
\usepackage{hyperref}
\newcommand{\prob}{\mathbb P}

\begin{document}

\title{First exit time for a discrete time parallel queue}

\author{Zbigniew Palmowski
}

\institute{Department of Applied Mathematics, Wroclaw University of Science and Technology, Wroclaw, Poland. \email{zbigniew.palmowski@pwr.edu.pl}
       }

\date{\today}

\maketitle

%
\section{Introduction}
We consider a discrete time parallel queue, which is a two-queue network, with
batch arrivals and services. At the beginning of the $n$th time-slot, $A_n$ customers
arrive to both queues. Then after this,
$S_n^i\in \mathbb{N}\cup\{0\}$ customers can be served in slot $n$ by server $i$
for $i=1,2$. We assume that $\{A_n\}_{\{n\in \mathbb{N}\cup\{0\}\}}$, $\{S_n^i\}_{\{n\in \mathbb{N}\cup\{0\}\}}$
are independent sequences of i.i.d. random variables with support in $\mathbb{N}\cup\{0\}$ and $\mathbb{E}A_n<\mathbb{E}S_n^i$
for $i=1,2$.
Let $Q_n^i$ with $i=1,2$ be the queue length after the service $S_{n-1}^i$ and before the arrival $A_n$.
Then it satisfies the Lindley recursion
$Q_{n+1}^i=(Q_n^i+A_n-S_n^i)_+$ for $i=1,2$
and its stationary law $(Q_\infty^1, Q_\infty^2)$ is given by
$(\max_{n\in \mathbb{N}\cup\{0\}}T_n^1, \max_{n\in \mathbb{N}\cup\{0\}} T_n^2)$ where
$T_n^i=\sum_{k=1}^{n} (A_k-S_k^i),\quad T_0^i=0, \quad i=1,2.$
For
\[\mathcal{A}_n=\sum_{k=1}^n A_k, \quad \mathcal{S}_n^i=\sum_{k=1}^n S_k^i, \quad i=1,2,\]
and $x,y\in \mathbb{N}\cup\{0\}$, we are interested in
\begin{eqnarray*}
H(x,y)&=&\prob (Q_\infty^1>x, Q_\infty^2>y)\\
&=&\prob(\text{there exists $n\in \mathbb{N}$ such that}\quad
\mathcal{A}_n>\max\{x+\mathcal{S}_n^1, y+\mathcal{S}_n^2\}).
\end{eqnarray*}
Note that $H(x,y)$ equals the probability that the first entrance time of the two-dimensional random walk $\{(\mathcal{A}_n-\mathcal{S}_n^1, \mathcal{A}_n-\mathcal{S}_n^2), n\in \mathbb{N}\cup\{0\}\}$ to the set $[x, \infty]\times [y, \infty]$ is finite.
In the context of actuarial science, $H(x,y)$ corresponds to ruin probability for the
a discrete-time two-dimensional insurance risk model where each business line faces
common simultaneous losses $A_k$ and $S_k^i$, for $i=1,2$, is the premium derived within the $k$th period of time
(see \cite{Serban} for similar considerations).

This two-queue network is a very special case of a discrete time network, but we hope that it provides very interesting insight
in understanding the more general picture (see [1,2]). Very natural generalizations can concern larger amount of parallel queues,
the case where apart from the common input each queue has its own batch arrivals etc.

\section{Problem statement}
The natural questions concern the exact form of $H(x,y)$ and its asymptotics for large values of $x$ and $y$.\\
{\bf Exact expression for $H(x,y)$.}
We have $H(x,y)=\sum_{m>x, n>y}p(m,n)$ for the stationary distribution $p(m,n)=\prob(Q_\infty^1=m, Q^2_\infty=n)$.
Moreover, $p(m,n)$ is a stationary distribution of the discrete time Markov chain $\{(Q_n^1, Q^2_n), n\in \mathbb{N}\cup\{0\}\}$,
hence it satisfies a classical balance equation. This probably can be used to get $H(x,y)$ at least
in two cases.
In the first case, $A_n, S_n^i$, $i=1,2$, are $0-1$ random variables
with $\prob (A_n=1)<\min\{\prob (S_n^1=1), \prob (S_n^2=1)\}$.
In the second case, we can assume that $A_n, S_n^i$ ($i=1,2$) have geometric laws, that is,
\[\prob (A_n=k)=\alpha (1-\alpha)^{k}, \qquad \prob (S_n^i=k)=\beta_i (1-\beta_i)^{k}, \quad i=1,2,\quad k=0,1,2,3,\ldots\]
and $\frac{1-\alpha}{\alpha}<\min\left\{\frac{1-\beta_1}{\beta_1},\frac{1-\beta_2}{\beta_2}\right\}$.
\\
{\bf Light-tailed asymptotics.}\\
Under light-tailed assumptions we look for the asymptotics of
$H(n\eta_1 ,n\eta_2 )$ as $n\rightarrow \infty$, where
$\eta_i>0$, i=1,2.
\\
{\bf Heavy-tailed asymptotics.}
Here we assume the complementary condition that the arrival sizes $A_n$ are strongly subexponential, that is, that
\[\lim_{n\to\infty}\frac{\prob(A_1+A_2>n)}{\prob(A_1>n)}= 2\quad\text{and}\quad \sum_{k=0}^n\prob(A_1>n-k)\prob(A_1>k)\sim 2 \mathbb{E}A_1 \prob(A_1>n)\]
as $n\rightarrow +\infty$, where $f(n) \sim g(n)$ when $\lim_{n\to\infty}f(n)/g(n)=1$.
The goal is to use the principle of a single big jump to find (under some moment assumptions put on the service capacities $S_n^i$
served in one slot), 
the asymptotics of $H(n\eta_1 ,n\eta_2)$ as $n\to \infty$.
\section{Discussion}
{\bf Exact expression for $H(x,y)$.} We believe that one should start from the balance equation. For example
in the second case we have
\[p(m,n)=\sum_{k=-\infty}^m\sum_{l=-\infty}^n\sum_{l_1-s=k \atop l_2-s=l} p(m-k, n-l)\alpha(1-\alpha)^s\beta_1(1-\beta_1)^{l_1}\beta_2(1-\beta_2)^{l_2}. \]
Moreover, one can use the bivariate moment generating function
$\mathbf{H}(z,w)=$\linebreak $\sum_{m,n=0}^\infty z^mw^np(m,n)$ to simplify above equations and try to identify the solution (see \cite{Cohen} and  \cite{FayK} for similar considerations).
\\
{\bf Light-tailed asymptotics.}\\
The case when  $A_n, S_n^i$, $i=1,2$, are $0-1$ random variables was analysed in \cite{Faybook}. For the case
of an arbitrary direction, \cite{BorMoga} derived similar but finer asymptotics.
However, it seems that they are not fully proved. It would be nice if the conjecture can be fully
proved. For this, it might be necessary to restrict the direction for the tail asymptotic
because the polynomial pre-factor may be different from $n^{-1/2}$ depending on the direction
as studied in \cite{Miyazawa2,Miyazawa1}.

To get exponential asymptotics for general distributions of $A_n, S_n^i$, $i=1,2$ (for example geometric ones)
one can follow
\cite{BorMog,Collamore,Miyazawa2,Miyazawa1}.
That is, for $\varphi \in \mathbb{R}^2$ we denote $\varphi(\vartheta)=\mathbb{E} \exp\{ <\mathbf{Q}_n, \vartheta>\}$.
We assume that there exists a solution $(\gamma, s) \in \mathbb{R}^2\times (0,\infty)$ of the Cram\'er equation:
$\varphi(\gamma)=1,$ $\varphi^\prime(\gamma)=\eta s,$
where $\eta=(\eta_1, \eta_2)$ and
$\varphi^\prime(\gamma)=\left(\frac{\partial \varphi (\gamma)}{\partial \gamma_1}, \frac{\partial \varphi (\gamma)}{\partial \gamma_2}\right)$.
Then one can conjecture that (under some additional assumptions)
\begin{equation*}
H(n\eta_1 ,n\eta_2)\sim C n^{-1/2}e^{-<\gamma, \eta>n}
\end{equation*}
for some constant $C$. One can think of
identifying the constant $C$ for some particular laws of arrival and service sizes as well.\\
{\bf Heavy-tailed asymptotics.}
Assume that $S_n^i\geq 1$, $i=1,2$. We believe that one can use an idea given in \cite{FKPR} representing $H(x,y)$
as a crossing probability of an increasing and random barrier by a random walk $\mathcal{A}_n$ and then apply \cite{FPZ} to prove that
\begin{equation*}
H(n\eta_1 ,n\eta_2)\sim \sum_{k=0}^\infty \prob(A_1>\max\{ n\eta_1+k\mathbb{E} S_1, n\eta_2+k\mathbb{E} S_2 \})\quad\text{as $n\rightarrow+\infty$.}
\end{equation*}

\begin{acknowledgements}
The research of Zbigniew Palmowski is partially supported by Polish National Science Centre Grant
No. 2018/29/B/ST1/00756 (2019-2022).
\end{acknowledgements}

\bibliographystyle{abbrv}
	
\end{document}